\documentclass[11pt]{amsart}
\usepackage[normalem]{ulem}
\usepackage[usenames]{color}
\usepackage{graphicx,amscd,psfrag,amssymb}
\usepackage{citesort}

\DeclareMathOperator{\Inv}{Inv} \newcommand{\R}{\mathbb R}
 \newcommand{\Z}{\mathbb Z}
\newcommand{\bs}{\backslash}\newcommand{\ep}{\varepsilon}

\DeclareMathOperator{\Int}{Int}
\DeclareMathOperator{\cl}{cl}

\newcommand{\p}{\mathcal{P}}

\newcommand{\Con}{\operatorname{Con}}
\newcommand{\Id}{\operatorname{Id}}

\theoremstyle{plain} \newtheorem{thm}{Theorem}
 \newtheorem{prop}[thm]{Proposition}
\newtheorem{lemma}[thm]{Lemma}

\theoremstyle{definition} \newtheorem{defn}[thm]{Definition}

\newtheorem{ex}[thm]{Example} 

\theoremstyle{remark}

\begin{document}

\title{Symbolic dynamics for nonhyperbolic systems}

\author{David Richeson}   \author{Jim Wiseman} \address{Dickinson
College\\ Carlisle, PA 17013} \email{richesod@dickinson.edu} 
\address{Agnes Scott College \\ Decatur, GA 30030}
\email{jwiseman@agnesscott.edu}


\begin{abstract}
We introduce index systems, a tool for studying isolated invariant sets of dynamical systems that are not necessarily hyperbolic.  The mapping of the index systems mimics the expansion and contraction of hyperbolic maps on the tangent space, and they may be used like Markov partitions to generate symbolic dynamics. Every continuous dynamical system satisfying a weak form of expansiveness possesses an index system. Because of their topological robustness, they can be used to obtain rigorous results from computer approximations of a dynamical system.
\end{abstract}

\maketitle

\section{Introduction}

Hyperbolicity is one of the most important ideas in dynamical systems.  Every hyperbolic diffeomorphism admits a Markov partition---a finite collection of rectangles that stretch or shrink in different directions, and map nicely onto each other. The interaction of rectangles under a single application of the map is sufficient to generate symbolic dynamics, which in turn gives global information about the dynamical system.  

In this paper we introduce index systems, a topological generalization of hyperbolicity and Markov partitions. Index systems are composed of finitely many index pairs, a fundamental object from Conley index theory. The mapping of the index pairs mimics the expansion and contraction of the Markov rectangles, and they may be used to generate symbolic dynamics. 

The benefit of this topological approach is that it applies in much more general situations than do Markov partitions. Hyperbolicity is a strong condition (requiring at least a manifold and a  differentiable map) that can be difficult to verify. However, index systems can be constructed on metric spaces in which the map is not differentiable. In particular, every discrete dynamical system that satisfies a weak form of expansiveness (one of the properties of every hyperbolic system) is guaranteed to have a nontrivial index system. Moreover, unlike Markov partitions, index systems are robust under slight perturbations of the map. This will enable us (in future work) to obtain rigorous results about dynamical systems from computer approximations. As we will see, however, the expense of this generality is that whereas a Markov partition produces a subshift of finite type, an index system generates a cocyclic subshift.

The paper is organized as follows.  In Section~\ref{sec:background} we recall necessary information about expansiveness and the Conley index.  We discuss index systems and their properties in Section~\ref{sec:systems}. We show how to use them to detect orbits of the dynamical system in Section~\ref{sec:detecting} and how they can be used to generate symbolic dynamics in Section~\ref{sec:symb}.  In Section~\ref{sec:existence} we show that all maps with a weak form of expansiveness have index systems. In Section~\ref{sec:examples} we give examples.

\section{Background}\label{sec:background}

Unless otherwise specified, throughout this paper we let every space $X$ be a compact metric space and every dynamical system $f:X\to X$ be a continuous map.  An \emph{orbit} of $f$ is a bi-infinite sequence $(x_{i}\in X:i\in\Z)$ with the property that $f(x_{i})=x_{i+1}$ for all $i\in\Z$. 

\subsection{Expansiveness}
\label{ssect:expansiveness}
A homeomorphism $f:X\to X$ is \emph{expansive} if there exists $\rho>0$ such that for any distinct points $x,y\in X$, there is an integer $n$ with $d(f^n(x),f^n(y))>\rho$. In other words, if the orbits of two points stay close together for all time, then they must be the same point.  Expansiveness is a strong form of sensitive dependence on initial conditions since any two distinct points must eventually move apart in either forward or backward time. In practical terms, this means that any small initial measurement error will lead to large errors in predicting behavior.  It is not difficult to show that if ${S}$ is a hyperbolic invariant set for a diffeomorphism $f$, then $f$ restricted to ${S}$ is expansive (\cite{KH,Ro}).

Although expansiveness is defined as a metric property, on compact spaces it is independent of the metric (of those compatible with the topology), and there is a very simple and useful topological characterization. Moreover, the topological definition can be extended trivially to continuous maps. To state it we need some preliminary definitions.

We begin with the definition of an {isolated invariant set}, a notion that is also central to Conley index theory and will be discussed in more detail in Section~\ref{sec:conley}.  A set $S\subset X$ is an \emph{isolated invariant set} for a continuous map $f:X\to X$ provided there is a compact set $I$ with $S=\Inv I\subset\Int I$ (where $\Inv I$ denotes the maximal invariant subset of $I$).  The set $I$ is called an \emph{isolating neighborhood} for $S$. Examples of isolated invariant sets include hyperbolic periodic
orbits, attractors, and the invariant Cantor set inside the Smale horseshoe.

For any dynamical system $f:X\to X$, let $f\times f:X\times X\to X\times X$ be the map $(f\times f)(x_{1},x_{2})=(f(x_{1}),f(x_{2}))$. Let $1_{X}=\{(x,x)\in X\times X\}$ denote the diagonal of $X\times X$.

\begin{defn}
A continuous map $f:X\to X$ is {\em expansive} if the diagonal $1_{X}\subset X\times X$ is an isolated invariant set with respect to $f\times f$. 
\end{defn}

Note that when $f$ is a homeomorphism this definition of expansive is equivalent to the metric definition (\cite[Def.~11.4]{A}).

We have defined expansiveness for a map $f:X\to X$.  
If ${S}$ is an isolated invariant subset of $X$ we say that {\em $f$ is expansive on ${S}$} if $f$ restricted to ${S}$ is expansive. Note that this is equivalent to the condition that the set $1_{S}=\{(x,x) : x\in {S}\}\subset X\times X$ is an isolated invariant set for $f\times f:X\times X\to X\times X$ (\cite[Ex.~11.5]{A}).

\subsection{Conley index}
\label{sec:conley}

The discrete Conley index is a powerful topological tool for studying
isolated invariant sets.  Roughly speaking, the Conley index assigns to each isolated invariant
set for $f$ a pointed topological space and a base-point preserving map on it,
$f_P$, which is unique up to an equivalence relation.  By studying
the simpler map $f_P$ we can draw conclusions about the original map
$f$.  Our discussion of the discrete Conley index is based on that in
\cite{FR}, where one can find more details and proofs of the theorems
below.

Ideally one would like to place an isolated invariant set, $S$,
between levels of a filtration.  In other words, we would like sets
$N_{0}\subset N_{1}$, both of which map into their interiors, such that
$S=\Inv(N_{1}\bs N_{0})$.  In practice this may not be possible.
Instead we must settle for a topological pair that behaves locally like
a filtration.

\begin{defn}\label{filpair}
Let $S$ be an isolated invariant set and suppose $L \subset N$ are
compact sets.  The pair $(N,L)$ is an {\em index pair} for
$S$ provided $N$ and $L$ are each the closures of their interiors and
\begin{enumerate}
\item
$\cl (N\bs L)$ is an isolating neighborhood for $S$,
\item
$L$ is a neighborhood of the {\em exit set}, $N^{-}=\{x\in N: f(x)
\notin \Int N\}$, in $N$, and
\item
$f(L) \cap \cl (N\bs L) = \emptyset$.
\end{enumerate}
\end{defn}

\begin{figure}[ht]
\centering
\psfrag{L}{$L$}
\psfrag{fL}{$f(L)$}
\psfrag{fN}{$f(N)$}
\psfrag{N}{$N$}
\psfrag{Lp}{$L^{\prime}$}
\psfrag{fLp}{$f(L^{\prime})$}
\psfrag{fNp}{$f(N^{\prime})$}
\psfrag{Np}{$N^{\prime}$}
\includegraphics{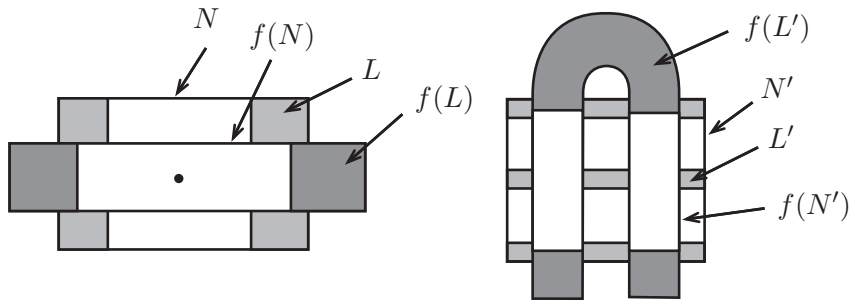}
\caption{Index pairs: $(N,L)$ for a fixed point saddle,
$(N^{\prime},L^{\prime})$ for the horseshoe Cantor set, and their
images
under $f$.}
\label{fig:conley}
\end{figure}

The definition of index pair given above was introduced in
\cite{FR} (where they were called filtration pairs).  This definition is similar to those in \cite{BF,E,Mr1,Sz,RS}.  We see examples
of index pairs in Figure \ref{fig:conley}.

\begin{figure}[ht]
\centering
\psfrag{NL}{$N_{L}$}
\psfrag{NpLp}{$N^{\prime}_{L^{\prime}}$}
\psfrag{[L]}{$[L]$}
\psfrag{[Lp]}{$[L^{\prime}]$}
\includegraphics{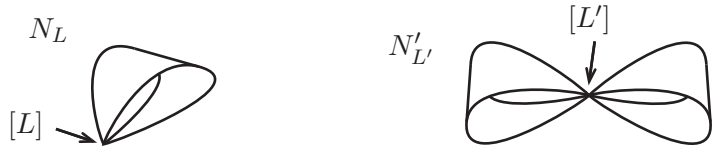}
\caption{The pointed spaces $N_{L}$ and $N^{\prime}_{L^{\prime}}$
obtained from the filtration pairs for the saddle and the horseshoe.}
\label{fig:conley2}
\end{figure}

Given a neighborhood $U$ of an isolated invariant set $S$ there exists
an index pair $P=(N,L)$ with $N\bs L\subset U$.  Given such an
index pair we form the pointed space $N_{L}$ by collapsing $L$ to
a point, $[L]$ (see Figure \ref{fig:conley2}).  Thus $f$ induces a
continuous map $f_{P}:N_{L}\to
N_{L}$.  The base point $[L]$ is an attracting fixed point of $f_{P}$.  The choice of index pairs for $S$ is not
unique,
and different choices can lead to topologically different maps $f_P$.
However, the choice is unique up to  \emph{shift
equivalence} (we will use $\cong$ to denote shift equivalence).  The resulting equivalence class is called the {\em
homotopy Conley index}.  For a definition of shift equivalence and a
proof of the facts in this paragraph see \cite{FR}.

In practice, one may wish to compute the relative homology (or
cohomology) of the pointed space and
the induced map $(f_{P})_{*}:H_{*}(N_{L},[L])\to H_{*}(N_{L},[L])$.
Again, this map is unique up to shift equivalence.  We call this equivalence class the \emph{homology Conley index}, $\Con_{*}(S)$. (Note: we may apply the Leray functor to $(f_{P})_{*}$ to obtain an automorphism of a graded group, which is an invariant for $S$ (see \cite{Mr1}).)

For instance, the fixed point saddle in Figure \ref{fig:conley} has
$\Con_{*}(S)\cong (\Z,\Id)$ where the $\Z$ is in dimension 1.  The invariant Cantor set in the
horseshoe has $\Con_{*}(S)\cong 0$ (from a homological perspective the
two arms of the horseshoe cancel one another). Note that that if $\Con_{*}(S)\not\cong 0$, then $S\ne\emptyset$ (this is called Wa\.{z}ewski's theorem), but the horseshoe shows that the converse is not true.

An important feature of the  Conley index is that if $I$ is an isolating neighborhood for $f$, then $I$ is also an isolating neighborhood for $\Inv(I,g)$ for any $g$ that is  $C^0$-close to $f$, and in this case $\Inv(I,f)$ and $\Inv(I,g)$ have the same Conley index (\cite{FR}).  
This robustness allows one to compute the Conley index from a suitable numerical approximation of the map (\cite{Ka,Mi,Mr2,Mr3}).

\section{Index systems}\label{sec:systems}

For differentiable maps it is frequently useful to look at the behavior of orbits on the manifold and also the behavior of the derivative on the tangent space following these orbits.  
Indeed, this is precisely what one does with hyperbolicity. As we see in this section, for an expansive map $f:X\to X$ we may use $f\times f:X\times X\to X\times X$ to obtain a nice analogue of the differentiable situation. We will see that slices of a neighborhood $N$ of the diagonal (that is,  intersections of $\{x\}\times X$ with $N$) are the analogues of the tangent spaces and it is from these slices that we build the index system.

Suppose ${\mathcal I}=(I_i:i\in\Z)$ is a sequence of compact sets in $X$ (in practice $\{I_{i}:i\in \Z\}$ will be a finite set of sets). We say that an orbit $(x_{i})$ \emph{follows} ${\mathcal I}$ if $x_{i}\in I_{i}$ for all $i\in \Z$.  We call such a sequence, ${\mathcal I}$, an \emph{isolating neighborhood chain} if for any orbit $(x_{i})$ that follows ${\mathcal I}$, $x_{i}\in\Int(I_i)$ for all $i\in\Z$.  For ${\mathcal I}$ an isolating neighborhood chain,  let $\Inv_0({\mathcal I}) \subset \Int(I_0)$ denote the set of points whose orbits follow ${\mathcal I}$; that is, \[\Inv_{0}({\mathcal I})= \{x\in I_0 : \text{there exists an orbit $(x_i)$, with $x_0=x$, that follows $\mathcal I$}\}.\]  Note that ${\Inv_0}({\mathcal I})$ is not, in general, $f$-invariant.

\begin{figure}[ht]
\centering
\psfrag{N1}{$(N_{a_i}, L_{a_i})$}
\psfrag{N2}{}
\psfrag{N3}{$(N_{a_{i+1}}, L_{a_{i+1}})$}
\psfrag{N4}{}
\psfrag{N5}{$(N_{a_{i+2}}, L_{a_{i+2}})$}
\psfrag{N6}{}
\includegraphics{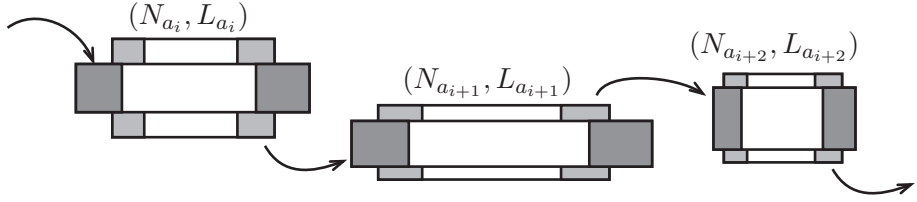}
\caption{$(N_{a_i}, L_{a_i})$ precedes $(N_{a_{i+1}}, L_{a_{i+1}})$ which precedes $(N_{a_{i+2}}, L_{a_{i+2}})$.} 
\label{fig:fsystems}
\end{figure}

An index system is a collection of compact pairs.  Each is similar to an index pair, but instead of necessarily mapping to itself under $f$, it maps to one or more of the pairs in the index system (see Figure \ref{fig:fsystems}).  More precisely, we have the following definition.

\begin{defn}
An {\em index system} is a finite collection of compact pairs, $\p=\{P_a=(N_a,L_a):a\in{\mathcal A}\}$, such that 
\begin{enumerate}
\item
for each $a\in\mathcal A$, there exists at least one $b\in\mathcal A$ such that $P_a$ \emph{precedes} $P_b$, that is, such that
	\begin{enumerate}
	\item $L_a$ is a neighborhood of the exit set, $N_{ab}^- = \{x\in N_a: f(x) \not\in \Int N_b\}$, in $N$, and
	\item $f(L_a)\cap \cl(N_b\bs L_b)=\emptyset$, and
	\end{enumerate}
\item
any sequence $(I_i=\cl(N_{a_i}\bs L_{a_i}): P_{a_i}\text{ precedes }P_{a_{i+1}}\text{ for all }i\in\Z)$ is an isolating neighborhood chain.
\end{enumerate}
\end{defn}

For each $a$, we form the pointed space ${N_a}_{L_a}$ by collapsing $L_a$ to
a point, $[L_a]$.  If $P_a$ precedes $P_b$, then $f$ induces a continuous map $f_{a,b}:{N_a}_{L_a} \to {N_b}_{L_b}$.

We can think of an index system as a directed graph.  Each pointed space ${N_a}_{L_a}$ is a vertex, with an edge from ${N_a}_{L_a}$ to ${N_b}_{L_b}$ if $P_a$ precedes $P_b$.  (The induced maps $f_{a,b}$ are analogues of the induced maps on the tangent bundle for differentiable systems.) We say that a finite or infinite sequence $(a_i)$ is {\em allowable} if it corresponds to a path in the graph, that is, if $P_{a_i}$ precedes $P_{a_{i+1}}$ for all $i$.

Given an index system $\p$, we define its {\em invariant set} by $\Inv(\p)=\bigcup \Inv_0((I_i=\cl(N_{a_i}\bs L_{a_i})))$, where the union is over all allowable sequences $(a_i:i\in\Z)$.  Thus $\Inv(\p)$ is the set of points on orbits following an allowable sequence of pairs. Let $\p$ be an index system.  Define $\Inv^m(\p)$ to be the set of $x$ such that there exists an orbit segment $(x_i)_{i=-m}^m$ with $x_0=x$ and $x_i \in \cl(N_{a_{i}}\bs L_{a_{i}})$ for some finite allowable sequence $(a_i)_{i=-m}^m$.

\begin{lemma}\label{lem:inv}
Let $\p$ be an index system.  Then $\Inv(\p) = \bigcap_{m=0}^\infty\Inv^m(\p)$.
\end{lemma}

\begin{proof}
The proof is based on that of \cite[Prop.~2.2]{FR}, which is the same result for index pairs instead of index systems.  It is obvious that $\Inv(\p) \subset \bigcap_{m=0}^\infty\Inv^m(\p)$.  To prove the opposite inclusion, let $x$ be an element of $\bigcap_{m=0}^\infty\Inv^m(\p)$; we must show that $x \in \Inv(\p)$.

For an allowable sequence $(a_{-1}, a_0)$ with $x \in \cl(N_{a_0}\bs L_{a_0})$, define $X_1^{(a_{-1}, a_0)}=f^{-1}(x) \cap \cl(N_{a_{-1}}\bs L_{a_{-1}})$.  Then, for $(a_{-k},\dots, a_0)$ an allowable sequence with $x \in \cl(N_{a_0}\bs L_{a_0})$, inductively define $X_k^{(a_{-k}, \dots,a_0)}=f^{-1}(X_{k-1}^{(a_{-k+1}, \dots,a_0)}) \cap \cl(N_{a_{-k}}\bs L_{a_{-k}})$.  Each $X_k^{(a_{-k}, \dots,a_0)}$ is compact and $f(X_k^{(a_{-k}, \dots,a_0)}) \subset X_{k-1}^{(a_{-k+1}, \dots,a_0)}$.  Set $X_k = \bigcup X_k^{(a_{-k}, \dots,a_0)}$, where the union is over all allowable sequences $(a_{-k},\dots, a_0)$ with $x \in \cl(N_{a_0}\bs L_{a_0})$; then $X_k$ is a non-empty, compact set with $f(X_k) \subset X_{k-1}$.  

Now define $Y_k = \bigcap_{n\ge1}f^{n+1}(X_{n+k})$.  As the intersection of a nested sequence of non-empty compact subsets of $X_k$, $Y_k$ is a compact, non-empty subset of $X_k$, so $f(Y_1)=\{x\}$ and for $k>1$,
\[ f(Y_k)=\bigcap_{n\ge1}f^{n+1}(X_{n+k}) = \bigcap_{n\ge2}f^{n}(X_{n+k-1}) = Y_{k-1}. \]
Finally, define $x_{-1}$ to be any point of $Y_1$, and inductively define $x_{-k}$ to be any point in $Y_k$ with $f(x_{-k})=x_{-k+1}$.  If we set $x_k=f^k(x)$ for $k\ge0$, then $(x_i)_{i=-\infty}^\infty$ is an orbit following an allowable sequence, so $x\in \Inv(\p)$.

\end{proof}

The following result follows immediately from Lemma \ref{lem:inv}.

\begin{prop}
If $\p$ is an index system, then $\Inv(\p)$ is a compact invariant set.
\end{prop}

\begin{defn} 
Let ${S}\subset X$ be an isolated invariant set.  An index system $\p$ is an {\em index system for ${S}$} if
$\Inv(\p) = {S}$.
\end{defn}

In Section~\ref{sec:existence} we present situations in which we can guarantee the existence of index systems, and we give a procedure for constructing them that could be implemented on a computer.

\section{Detecting Orbits}
\label{sec:detecting}

In this section we show how to prove the existence of an orbit following a given allowable sequence of sets $(\cl(N_{a_i}\bs L_{a_i}))$. If $f$ had the shadowing property and the sets $\cl(N_{a_i}\bs L_{a_i})$ were small enough, then we could guarantee the existence of an orbit that follows this sequence (\cite[\S9.3]{Ro}).  However, the shadowing property is difficult to verify in practice. Instead, we use the Conley index to verify that there is an $f$-orbit that ``shadows'' the sequence in the sense that each iterate is in the appropriate set $\cl(N_{a_i}\bs L_{a_i})$.

An allowable sequence $(a_i)$ yields a directed system 
\[\begin{CD} \cdots \longrightarrow H_{*}(P_{a_i}) @>f_{a_{i},a_{i+1}*}>>H_{*}(P_{a_{i+1}}) @>f_{a_{i+1},a_{i+2}*}>>H_{*}(P_{a_{i+2}}) \longrightarrow \cdots
\end{CD}.\] 
The sequence $(P_{a_i})$ has a \emph{nonzero orbital Conley index} if any finite composition $f_{a_{n-1},a_{n}*}\circ \dots \circ f_{a_{m+1},a_{m+2}*} \circ f_{a_{m},a_{m+1}*}$ is nonzero (note that if we take coefficients in a field and $H_{*}(P_a)$ is finitely generated for each $a$, then the maps are simply matrices).

We have the following theorem, which is the orbital analogue of Wa\.{z}ewski's theorem.

\begin{thm}\label{thm:waz}
Let $\p$ be an index system for the isolated invariant set ${S}$, and let $(a_i)$ be an allowable sequence. If $(P_{a_i})$ has a nonzero orbital Conley index, then there is an orbit in ${S}$ following $(\cl(N_{a_i}\bs L_{a_i}))$. If $f$ is expansive and the diameters of the sets $\cl(N_{a}\bs L_{a})$ are sufficiently small for all $a$, then this orbit is unique.
\end{thm}

\begin{proof}
If the set $\bigcap_{i=-m}^m f^{-i} (\cl(N_{a_i}\bs L_{a_i}))$ were empty, then the induced map would send every point to the basepoint, and the corresponding map on homology would be zero.  Since $(P_{a_i})$ has a nonzero orbital Conley index, the set $\bigcap_{i=-m}^m f^{-i}(\cl (N_{a_i}\bs L_{a_i}))$ is nonempty for all $m$.  The same argument used to prove Lemma~\ref{lem:inv} then shows that there is an orbit in ${S}$ following $(\cl(N_{a_i}\bs L_{a_i}))$.

Suppose $f$ is expansive on ${S}$ and $(x_{i})$ and $(y_{i})$ are two orbits following $(\cl(N_{a_i}\bs L_{a_i}))$. Because $f$ is expansive, $1_{{S}}$ is an isolated invariant set for $f\times f$ with some isolating neighborhood $I$. Assuming the diameters of $\cl(N_{a_i}\bs L_{a_i})$ are small enough, that means that $(f\times f)$-orbit $((x_{i},y_{i}))$ remains in $I$ for all $i$. But this implies that for all $i$, $(x_{i},y_{i})\in 1_{{S}}$, or equivalently $x_{i}=y_{i}$.
\end{proof}

\section{Symbolic dynamics}
\label{sec:symb}

Many authors have used the Conley index or related techniques (for example, Easton's windows (\cite{Eas81})) to detect symbolic dynamics; an incomplete list of references is \cite{DFT,Gi,Mr2,Wis1,WZ,Z,CKM,WisSquare,Sz96,MM,Eas81,GR,ZG,vVW,Gia,Eas75,Gid99}. Most treatments resemble index systems, except that all of the index pairs are pairwise disjoint. This is an important difference---disjoint index pairs make the conclusions stronger, but they are significantly less useful because they cover less of the space. In \cite{Gid99}, Gidea developed a notion similar to index systems, of an orbital Conley index for compact sets traveling within a prescribed sequence of neighborhoods, which he applied to detect periodic orbits and symbolic dynamics.  An important difference is that index systems are used to study invariant sets, while Gidea's orbital Conley index applies to non-invariant sets that travel along a bi-infinite sequence of boxes.  In addition, there is no result similar to Theorem~\ref{thm:systemsexist} for constructing the neighborhood sequence for Gidea's orbital Conley index, and the method for detecting symbolic dynamics requires that the neighborhoods be pairwise disjoint.  As we see in this section, we can relax this requirement with index systems.

The mapping of rectangles in a Markov partition generates a directed graph, and paths through this graph become the elements of a subshift of finite type. The elements of this subshift can be paired uniquely with points in the dynamical system. The situation for index systems is similar, but more complicated. Intuitively, we would like allowable sequences from our directed graph of pointed spaces to generate symbolic dynamics. However, we must use the Conley index (as in Section~\ref{sec:detecting}) to determine when there is a true orbit corresponding to this allowable sequence, and then we must show that each point corresponds to only one symbol sequence.

In order to put this in proper context we must introduce the notion of a cocyclic subshift, a generalization of sofic shifts and subshifts of finite type. See \cite{Kwa} for more information about cocylic subshifts and \cite{CKM,Wis1,WisSquare} for examples of cocyclic subshifts in other dynamical and Conley index contexts.

\begin{defn}
Let $G$ be a directed graph with vertices $\{1,\ldots,n\}$ and at most one edge from any vertex to any other.  Assign to each vertex $i$ a vector space $V_i$ and to the edge from $i$ to $j$ a linear transformation $T_{i,j}:V_i\to V_j$.  The {\em cocyclic subshift} associated to $G$ consists of all words $(\ldots, \omega_{i-1}, \omega_i, \omega_{i+1},\ldots)$ such that
\begin{enumerate}
	\item there is an edge from vertex $\omega_{i}$ to $\omega_{i+1}$ for all $i$, and
	\item any finite composition $T_{i+m-1,i+m}\circ \cdots \circ T_{i,i+1}$ is nonzero,
\end{enumerate}
together with the usual shift map.
\end{defn}

An index system generates the cocyclic subshift corresponding to the graph with vertex set $\mathcal A$, vector spaces $H_*(P_a)$, edges from $a$ to $b$ if $P_a$ precedes $P_b$, and linear transformations $f_{a,b*}$. In other words, we take the homology of our directed graph of pointed spaces. At each vertex we place the homology of the given pair,  $H_*(P_{a}) = H_*({N_a}_{L_a},[L_a])$, and we label each directed edge from $H_*(P_a)$ to $H_*(P_b)$ with the corresponding induced map on homology $f_{a,b*}$.  Then elements of the cocyclic subshift correspond to allowable sequences with nonzero orbital Conley index.

A further complication is that in general we do not get a conjugacy (or semiconjugacy) from $f$ to the cocyclic subshift because the sets $\{\cl(N_a\bs L_a)\}$ are usually not pairwise disjoint.  For example, say
that the period-two word $(1,2,1,2,\dots)$ is in the cocyclic
subshift, which implies that there is a point $x\in X$ such that $x\in
\cl(N_{1}\bs L_{1})$, $f(x)\in \cl(N_{2}\bs L_{2})$, $f^2(x)\in \cl(N_{1}\bs L_{1})$, and so on.  But if the two sets are not disjoint, then $x$
could be a fixed point.  The $\p$-itinerary
of a point is not necessarily unique.

One way around this issue is simply to remove overlapping index
pairs.  That is, take a maximal subgraph of the cocyclic subshift
graph consisting of vertices for which all of the corresponding sets $
\cl(N_a\bs L_a)$ are pairwise disjoint. 
 Then, as desired, we will get a semi-conjugacy from an isolated
invariant set in $X$ to the cocyclic subshift induced by the
subgraph.  The disadvantage of this approach is that we are losing
information about $f$ when we throw away vertices.

Another approach is to look at powers $f^n$.  In this case, only some,
not all, of the index pairs must be disjoint in order to detect
positive-entropy symbolic dynamics.
For example, let $\omega=(\omega_0,\dots,\omega_n)$ be a word in the
cocyclic subshift, and define ${S}_\omega\subset X$ to be the subset
$\Inv(\{x\in X : f^i(x) \in \cl(N_{\omega_i}), 0\le i\le n
\} , f^{n+1})$.  If, say, the sets $\cl(N_{1}\bs L_{1})$ and $\cl(N_{3}\bs L_{3})$ are disjoint, then the sets ${S}_{(1,1,1)}$ and ${S}_{(1,2,3)}$
are disjoint (even if $\cl(N_{1}\bs L_{1})$, $\cl(N_{2}\bs L_{2})$, and $\cl(N_{3}\bs L_{3})$ are not pairwise disjoint), so we can define symbolic
dynamics and get a semiconjugacy for $f^3$.

Using either method for a generic index system, we get only a semiconjugacy from the invariant
set to the cocyclic subshift. Thus, while we can get a lower bound
for the topological entropy of $f$, we cannot conclude that a periodic
word for the subshift actually corresponds to a periodic point for $f
$. However,  if $f$ is expansive and the index pairs are sufficiently small, then by the uniqueness of orbits guaranteed by Theorem~\ref{thm:waz}, we do get a conjugacy.
 Even in some cases when we do not have expansiveness (in particular, if the homology is nontrivial on exactly one level), we can use a version of the Lefschetz fixed point theorem (\cite{Sz96}) instead of Theorem
\ref{thm:waz} to detect periodic points for $f$.  For the sake of brevity, we omit the details.

There are further details on techniques for generating symbolic dynamics from index systems in \cite{RWindexsymb}.

\section{Existence of index systems}
\label{sec:existence}

As we have said, hyperbolic diffeomorphisms admit Markov partitions. In fact, homeomorphisms that are expansive and have the shadowing property also admit Markov partitions (\cite{AH}), but expansivity alone is not sufficient. Furthermore, even when they exist, Markov partitions can be difficult to construct. 

One of the benefits of index systems is that under certain general circumstances their existence is guaranteed. Moreover, the proof is constructive---it requires only the ability to find a single index pair. There are computer algorithms that will do that (see \cite{DJM2}, for example).

It is immediate that every isolated invariant set ${S}$ has at least a trivial index system:  the system consisting of the single pair $(N,L)$, where $(N,L)$ is an index pair for ${S}$.  In the case that $f$ is expansive, we can do better.

\begin{thm}\label{thm:systemsexist}
Let $f$ be expansive on an isolated invariant set ${S}$.  Then there exists an index system $\p=\{(N_a,L_a)\}$ for ${S}$ of arbitrarily small diameter (that is, the diameter of $N_a$ is arbitrarily small for all $a$).
\end{thm}

\begin{proof}
The idea of the proof is simple.  Recall from Section~\ref{ssect:expansiveness} that since ${S}$ is expansive, the set $1_{S}=\{(x,x) : x\in {S}\} \subset X\times X$ is an isolated invariant set for $f\times f$.  Thus any neighborhood of $1_{S}$ contains an index pair $(N,L)$ for $1_{S}$ under $f\times f$.  By taking cross-sections, we get pairs $(N_x,L_x)=(\pi_2((\{x\}\times X) \cap N), \pi_2((\{x\}\times X) \cap L))$ of arbitrarily small diameter (where $\pi_2:X\times X\to X$ is the projection onto the second coordinate). We would like these pairs to make up the index system.  However, an index system must be finite, and in general this construction may give us infinitely many pairs, one for each $x\in{S}$.

By essentially discretizing the space and using the robustness of the Conley index we can modify the original index pair $(N,L)$ to yield only finitely many cross-sections.  Let $D$ be the metric on $X\times X$ given by $D((x,y),(x',y'))=\max(d(x,x'),d(y,y'))$.  By the definition of index pair, there exists an $\ep>0$ such that $D(1_{{S}},L)>\ep$, $D((f\times f)(L), \cl(N\bs L))>\ep$, and $D(N^-,\cl(N\bs L))>\ep$.  Pick $\delta<\ep/3$ such that for  any points $(x,y)$  and $(x',y')$ within $\ep$ of $N$, if $D((x,y),(x',y'))<\delta$, then $D((f\times f)(x,y),(f\times f)(x',y'))<\ep/3$.

Let $\{V_i\}$ be a finite collection of compact subsets of $X$ of diameter less than $\delta$ such that $\pi_1(N)\cup\pi_2(N)\subset\bigcup_i\Int(V_i)$.  This gives a product cover of $N$, $\{V_{ij}=V_i\times V_j\}$.  We now construct a new index pair for $1_{S}$, $(\tilde N,\tilde L)$, by setting $\tilde N = \bigcup_{\{ij:V_{ij}\cap N\ne \emptyset\}} V_{ij}$ and $\tilde L = \bigcup_{\{ij:V_{ij}\cap L\ne \emptyset\}} V_{ij}$.  (See Figure \ref{fig:discretized}.)  The set of cross-sections of $(\tilde N, \tilde L)$ is finite, and it is straightforward to check that it forms an index system for ${S}$.
\begin{figure}[ht]
\centering
\psfrag{X}{$X$}
\psfrag{1}{$1_{S}$}
\psfrag{N}{$(N_{a},L_{a})$}
\psfrag{p}{$x$}
\includegraphics{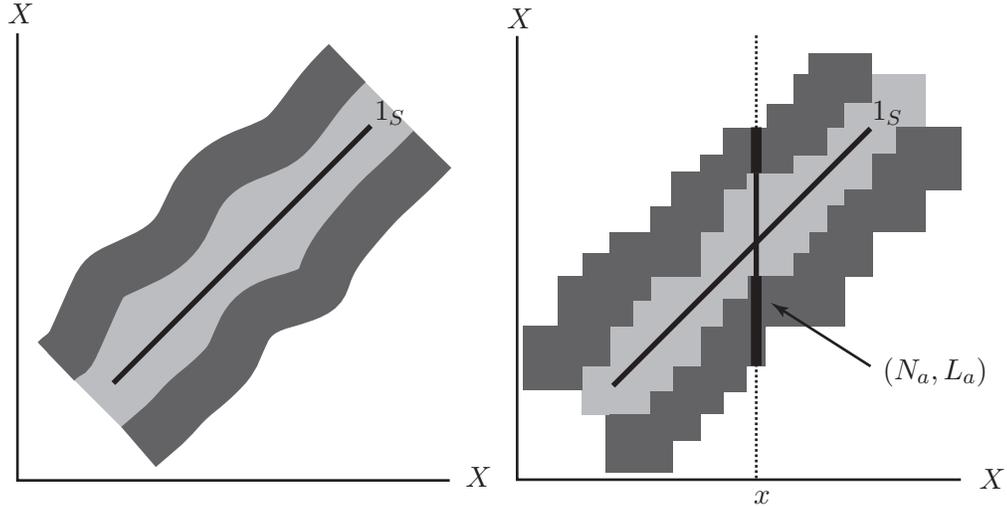}
\caption{The given index pair for $1_{{S}}$, a discretized index pair nearby, and one of the finite number of slices. }
\label{fig:discretized}
\end{figure}
\end{proof}

Theorem~\ref{thm:systemsexist} says that if $f$ is expansive on an isolated invariant set ${S}$, then it has an index system. Expansiveness is a strong condition that may be difficult to verify in practice. It turns out that a condition weaker than expansiveness  guarantees the existence of index systems. In Theorem~\ref{thm:systemsexist} we assumed that $1_{{S}}\subset X\times X$ was an isolated invariant set (i.e., that $f$ was expansive on ${S}$), but the technique for constructing the index system requires only that $1_{{S}}$ be contained in some isolated invariant set $\Lambda$. Then the slices of an index pair $(N,L)$ for $\Lambda$ give an index system for ${S}$.  At one extreme we have $\Lambda=1_{S}$ ($f$ is expansive), in which case we can find an index system with index pairs arbitrarily small. At the other extreme, for \emph{any} isolated invariant set ${S}$ we could take $\Lambda={S}\times{S}$; this would produce a trivial index system consisting of one index pair. There may be times when $1_{{S}}\subsetneq \Lambda\subsetneq {S}\times{S}$, and the resulting index system is useful.

If $1_{{S}}$ is not an isolated invariant set, then we will not be able to make the index pairs arbitrarily small, but that is not necessarily a big disadvantage. On the one hand, the smaller the sets $\cl(N_a\bs L_a)$, the stronger the conclusion of Theorem~\ref{thm:waz}, which is one reason the small diameters guaranteed by Theorem~\ref{thm:systemsexist} are important. Furthermore, the smaller the sets, the more disjoint pairs we have, which, as we have seen, make it easier to detect symbolic dynamics. On the other hand, small diameters can lead to a system with many small, individually unimportant pieces. This would produce a large graph, and thus a very complicated cocyclic subshift. So, in practice there is a trade-off involved in the size of the sets of the index system.

\section{Examples}\label{sec:examples}

\begin{ex}
Let $f:S^1\to S^1$ be the doubling map on the circle, considered as $\R/\Z$.  The collection $\p=\{(N_i,L_i)\}_{i=0}^9$ is an index system, where $N_i = [\frac{i-3-3\ep}{10},\frac{i+3+3\ep}{10}]$ and $L_i= [\frac{i-3-3\ep}{10},\frac{i-1-\ep}{10}] \cup [\frac{i+1+\ep}{10},\frac{i+3+3\ep}{10}]$.  $P_i$ precedes $P_j$ for $j=2i-1$, $2i$, or $2i+1\mod 10$ (see Figure \ref{fig:circle}).

\begin{figure}[ht]
\centering
\psfrag{0}{$P_{0}$}
\psfrag{1}{$P_{1}$}
\psfrag{2}{$P_{2}$}
\psfrag{3}{$P_{3}$}
\psfrag{4}{$P_{4}$}
\psfrag{5}{$P_{5}$}
\psfrag{6}{$P_{6}$}
\psfrag{7}{$P_{7}$}
\psfrag{8}{$P_{8}$}
\psfrag{9}{$P_{9}$}
\psfrag{fN0}{$(f(N_{0}),f(L_{0}))$}
\psfrag{N0}{$(N_{0},L_{0})$}
\psfrag{N1}{$(N_{1},L_{1})$}
\psfrag{N9}{$(N_{9},L_{9})$}
\psfrag{a}{$0$}
\psfrag{b}{$\frac{1}{10}$}
\psfrag{c}{$\frac{1}{5}$}
\psfrag{d}{$\frac{3}{10}$}
\psfrag{e}{$\frac{2}{5}$}
\psfrag{f}{$\frac{1}{2}$}
\psfrag{g}{$\frac{3}{5}$}
\psfrag{h}{$\frac{7}{10}$}
\psfrag{i}{$\frac{4}{5}$}
\psfrag{j}{$\frac{9}{10}$}
\includegraphics{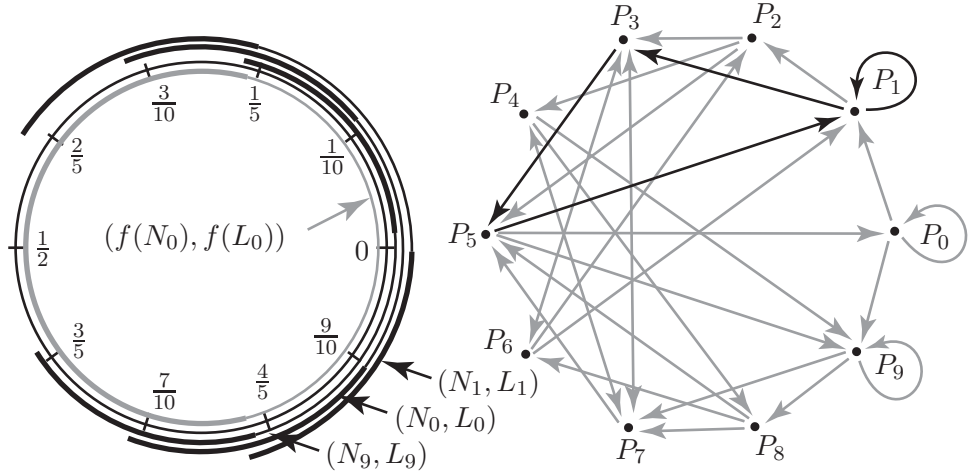}
\caption{An example showing that $P_{0}$ precedes $P_{9}$, $P_{0}$, and $P_{1}$, and the directed graph for this index system.}
\label{fig:circle}
\end{figure}

Each pointed space ${N_i}_{L_i}$ is homeomorphic to a circle. So, for every $i$, the only nonzero homology is $H_{1}(P_{i})=\Z$, and the induced maps are $f_{i,j*}=1$ if $P_i$ precedes $P_j$, and $f_{i,j*}=0$ otherwise.  Thus any concatenation of the maps $f_{(1,1,1)*}= f_{1,1*} \circ f_{1,1*} \circ f_{1,1*}$ and $f_{(1,3,5)*}= f_{5,1*} \circ f_{3,5*} \circ f_{1,3*}$ is nonzero.
Since the sets $\cl((N_1\bs L_1))$ and $\cl((N_5\bs L_5))$ are disjoint, so are the sets ${S}_{(1,1,1)}$ and ${S}_{(1,3,5)}$, and thus we see that the map $f^3:{S}_{(1,1,1)} \cup {S}_{(1,3,5)} \to {S}_{(1,1,1)} \cup {S}_{1,3,5}$ factors onto the full shift on two symbols.
\end{ex}

\begin{ex} \label{ex:tent}
Let $f:\R\to\R$ be the tent map given by $$f(x)=\begin{cases}  3x \text{, if } x\le 1/2,\\
3-3x \text{, if } x\ge 1/2.\end{cases}$$  
 
The map $f\times f$ has an index pair shown in Figure \ref{fig:product}.

\begin{figure}[ht]
\centering
\psfrag{a}{$\frac{1}{3}$}
\psfrag{b}{$\frac{2}{3}$}
\psfrag{c}{$1$}
\psfrag{N}{$(N,L)$}
\includegraphics{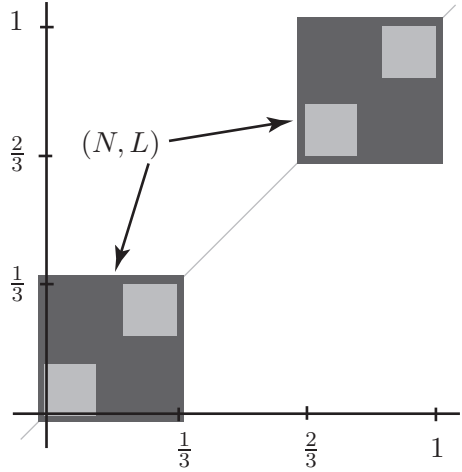}
\caption{An index pair for $f\times f$, where $f$ is the tent map.}
\label{fig:product}
\end{figure}

Taking slices we obtain the index system $\p=\{(N_i,L_i)\}_{i=1}^4$ for $f$ (as shown in Figure \ref{fig:tent}), where
$$\begin{array}{l}
N_1 = N_2 = [0-4\ep, \frac13+4\ep]\\
L_1 = [0-4\ep,0-\ep] \cup [\frac19 + \ep, \frac13 + 4\ep]\\
L_2 =  [0-4\ep,\frac29-\ep] \cup [\frac13 + \ep, \frac13 + 4\ep]\\
N_3 = N_4 = [\frac23 - 4\ep, 1 + 4\ep]\\
L_3 = [\frac23 - 4\ep, \frac23 -\ep] \cup [\frac79 + \ep, 1+4\ep]\\
L_4 = [\frac23 - 4\ep, \frac89 -\ep] \cup [1 + \ep, 1+4\ep].
\end{array}$$

Again, each pointed space ${N_i}_{L_i}$ is homeomorphic to a circle.  $P_1$ and $P_4$ precede $P_1$ and $P_2$, while $P_2$ and $P_3$ precede $P_3$ and $P_4$.  In homology, in dimension one, the induced maps are $f_{1,1*}= f_{1,2*} = f_{2,3*} = f_{2,4*} = 1$ and $f_{3,3*} = f_{3,4*} = f_{4,1*} = f_{4,2*} = -1$.  Thus the tent map restricted to $\Inv[0-\ep,1+\ep]$ factors onto the shift given by the graph in Figure \ref{fig:tent}, which is conjugate to the full shift on two symbols (\cite[\S2.4]{LM}).

\begin{figure}[ht]
\centering
\psfrag{N1}{$(N_{1},L_{1})$}
\psfrag{N2}{$(N_{2},L_{2})$}
\psfrag{N3}{$(N_{3},L_{3})$}
\psfrag{N4}{$(N_{4},L_{4})$}
\psfrag{P1}{$P_{1}$}
\psfrag{P2}{$P_{2}$}
\psfrag{P3}{$P_{3}$}
\psfrag{P4}{$P_{4}$}
\psfrag{a}{$\frac{1}{3}$}
\psfrag{b}{$\frac{2}{3}$}
\psfrag{c}{$1$}
\psfrag{d}{$\frac{1}{2}$}
\psfrag{e}{$\frac{3}{2}$}
\psfrag{f}{$1$}
\psfrag{g}{$-1$}
\includegraphics{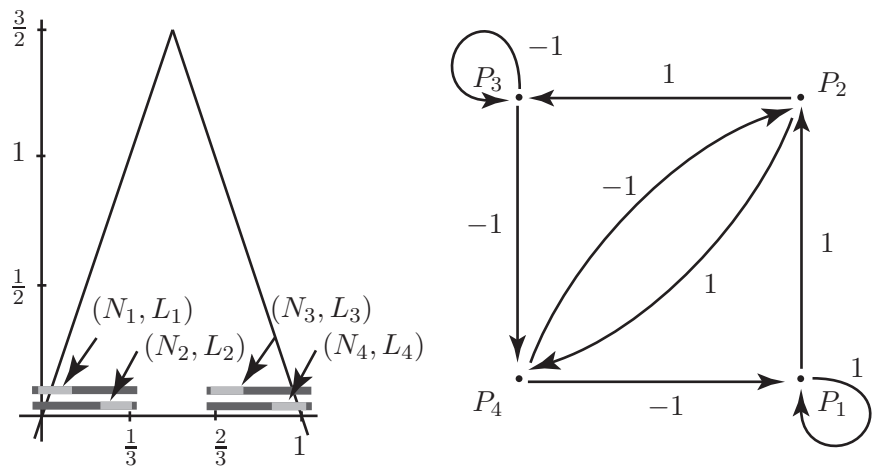}
\caption{An index system for the tent map, and the associated directed graph.}
\label{fig:tent}
\end{figure}
\end{ex}

\begin{ex}
Let $f$ be the tent map from Example~\ref{ex:tent}.  The pair $(N_0,L_0)=([-4\ep,1+4\ep],[-4\ep,-\ep]\cup[\frac{1}{3}+\ep,\frac{2}{3}-\ep]\cup[1+\ep,1+4\ep])$ is an index pair, and thus also a (trivial) index system.  The only nonzero homology is $H_1(N_0,L_0)=\Z^2$, and the induced map is given by $f_{00*}=\begin{pmatrix}1&-1\\1&-1\end{pmatrix}$.  Since $(f_{00*})^2$ is the zero matrix, the cocyclic subshift is empty.  Thus this index system fails to detect any invariant set.
\end{ex}

\bibliographystyle{amsplain} 
\bibliography{indexsystems} 
\end{document}